# A Relational Approach to Matroids, Simplicial Complexes and Finite Closures

Wayne E. Dick PhD.

Professor Emeritus, California State University at Long Beach

wayne.dick@csulb.edu, 1-562-256-4458

**Abstract**

The main result is Theorem [MAT 11]. If $\mu$ is a closure, then $\mathcal{M}(\mu)$ is a matroid with ground set $\mu$ and basis sets consisting of nonredundant covers of $\mu$, minimal subsets that determine $\mu$ using a closure algorithm from the theory of relational databases. The flat closures of matroid and simplicial complexes are compared to general closures using tools from Maier's algorithm to normalize databases [R1]. The flat closure operator for $\mathcal{M}(\mu)$ is called the span operator, $\mu^*$. This corresponds to linear span using the Extend by Closure recurrence.

## *[POV] Combining Points of View -*

The goal of this exposition is to merge two points of view regarding finite closures: The matroid / simplicial complex (MSC) approach [R2, R3, R4], and the relational database approach (RDB)[R1]. Both theories focus on functional relationships between rows and columns of tabular data. In MSC the matrix is the inspiration. In relational databases the tables of data that describe entities and their relationships motivates study. Both use the geometry and semantics of tabular data to infer structure on the models they represent.

Much of MSC theory seeks to model independence and dependence for values that can be computed by algebra. In the case of matroids the initial focus was linear algebra [R3]. One line of investigation for MSC seeks matrix representations of the closures over a field, ring or semiring [R2, R3, R4]. Matrices represent a matroid or simplicial complex if an independence preserving isomorphism can be established between an object's independent sets and the independent columns of the matrix. Many concepts from MSC, like flats and hereditary set systems, extend to the general closure operator. MSC theory is generally top-down. It moves from general structures like hereditary collections to specific dependencies.

Relational database theory approaches closures from a different direction. When data is stored in tables with named columns, natural functional dependencies occur between the columns because they are constructed to model relationships in real world systems. This is like the linear dependencies in matrices, but computable functional relationships are not stored in a database table. That would waste space. Computable values are manufactured by programs from independent data values when needed. Nonetheless certain functional relations exist between real world objects that must be reflected in the files that represent them. Consider the data generated for the columns named **{City, Year, Rainfall Total}.** There is a functional relationship between **{City, Year}** and **{Rainfall Total}**. On any given year a city receives so much rain. Nobody can compute this value with a formula, but nonetheless a functional relationship exists. Every time we encounter the city and year, (Los Angeles, 2020–2021 ) in a table row that includes its **Rainfall Total** , it must be paired with the value (5.00 inches), the quantity of rain fall that fell in Los Angeles in the year 2020-2021 . These existential functional relations determine closures that dictate which tables hold which columns of data. Relational database theory has classified these functional dependencies and their closures to protect data integrity and eliminate redundancy. There is a great deal of geometric and combinatorial structure to be learned from the study of these existential closures. In this context the functional dependency theory of RDB is a potential tool for exploring new properties of closures that relate directly to matroids and simplicial complexes. RDB theory is bottom-up. It proceeds from individual dependencies to general structures.

A shared aspect of closures with MSC is that they support greedy algorithms. The Extend by Closure recurrence given below is a greedy algorithm that builds the closure of any set using a "breadth first" approach. This works even though closures in general do not have an intrinsic rank function associated with their independent sets as is the case for matroids [R3, R4].

## *[CO] Closure Operators*

### Basic Concepts

Closures can be defined on general lattice structures. However, in this exposition only the subset lattice is considered.

**[CO 1]** Closure Operator: A *closure operator* is a function $\mu: 2^U \to 2^U$ given by $X \mapsto X\mu$ satisfying:

(C1) Inclusion: $X \subseteq X\mu$;

(C2) Monotonicity: $X \subseteq Y \Rightarrow X\mu \subseteq Y\mu$, and

(C3) Idempotence: $X\mu\mu = X\mu$.

Closures are unary operations and postfix notation is commonly used in this context. Thus $\mu(X) = X\mu$, like an exponent. This extends to normal function operations, like inverse image of a set of closed sets, $\{C_1, \ldots, C_k\}\mu^{-1}$.

**[CO 2]** Closed Sets : The fixed points of $\mu$ are *closed sets*.

**[CO 3]** Ground Set : The set $U$ is called the *ground* set of $\mu$. (This is taken from MSC)

## On Functions as Sets

A function $\alpha$ is a set of ordered pairs that satisfy $(X, Y) \,\&\, (X, Z) \in \alpha \Rightarrow Y = Z$. Set operations that apply to sets in general, union, intersection and set minus, also apply to functions. This can lead to confusion when we say $\alpha \subseteq \mu$ where both are functions. That will occur frequently in this paper. One item we prove here [**MAT 11**] is that every closure function $\mu$ is the ground set of an associated matroid $\mathcal{M}(\mu)$ with base sets that are minimal subsets (functions) in $\mu$ that generate $\mu$ using Extend by Closure [**EC 1**].

All the functions discussed here will map sets into sets. All objects are finite. The set $U$ and the closure $\mu$ will always be the universe that contains other sets and functions unless specifically noted. For, $\mu$ a closure, every subset $\alpha \subseteq \mu$ has the form $\alpha = \{(S, S\mu) \mid S \subseteq \text{dom}(\alpha)\}$, where $\text{dom}(\alpha) = U\alpha^{-1}$. We use the symbol $\subset$ to denote proper subset, and $\subseteq$ to denote the non-strict subset or equal.

## Common Structures Induced on the Ground Set

**[CO 4]** Key: A *key* $K$ is a minimal set regarding inclusion such that $(K, K\mu) \in \mu$. The set $K$ is a *key* for $C = K\mu$. The set of all keys of $\mu$ is denoted by $Ky\,\mu$.

**[CO 5]** Hereditary Collection: A collection of sets $H \subseteq 2^U$ is *hereditary* whenever $X \in H \ \& \ Y \subseteq X \Rightarrow Y \in H$.

**[CO 6]** *Theorem*: $Ky\ \mu$ is a hereditary collection for every closure $\mu: 2^U \to 2^U$.

Proof: Suppose there is nonempty $X \subset K$ and $X \notin Ky\ \mu$. We will show that $K \notin Ky\ \mu$.

$X \subset K$ and $X$ not a key, means, there are nonempty subsets $S, T, Y \subset K$ such that $X = S \cup T$, $K = S \cup T \cup Y$ and $S\mu = X\mu$. This is true because, $X$ is a proper subset of $K$ and $X$ is not a key.

Claim $(S \cup Y)\mu = K\mu$.

$S \cup Y \subseteq K \Rightarrow (S \cup Y)\mu \subseteq K\mu$.

To see the reverse inclusion $(K\mu \subseteq (S \cup Y)\mu)$, note the following:

$K\mu \subseteq (S\mu \cup T\mu \cup Y\mu)\mu$. $T\mu \subseteq X\mu = S\mu$. $\therefore K\mu \subseteq (S\mu \cup T\mu \cup Y\mu)\mu = (S\mu \cup Y\mu)\mu \subseteq (S \cup Y)\mu$.

Since a proper subset $S \cup Y$ of $K$ determines $K\mu$ the set $K$ is not a key. Thus if $K$ is a key every subset of $K$ must be a key. □

Note (a) The term key and independent set will be used interchangeably. A dependent set is a non-key.

Matroids and Simplicial Complexes start with hereditary collections of sets (independent) and induce special closed sets (flats) from their independent sets. This closure can be computed.

With general closers this process is reversed. The closure induces a hereditary collection of sets, keys. These independent sets are important, but they are not good to describe a dimension for general closures.

For example: In the database for a university, a student's data (full name, date of birth, and address) will uniquely identify the student and all the other personal data for the student at the university. It will also be a minimal set of data that does that. Hence, full name, date of birth and address form a key for the closure of the student's information. However, this key uses lots of characters, and the address is subject to change. Thus, the university will assign each student a unique ID number. The fact that one key has many fields and the other has one field tells little about the dimension of the closure space of student data.

The next theorem gives the relationship between a closed set $C$ and the rest of $U$.

**[CO 7]** *Theorem*: $C$ is closed regarding $\mu$ if and only if $p \in U \setminus C$ implies $C\mu \subset (C \cup \{p\})\mu$.

Proof: ($\Rightarrow$). $C$ closed means $C\mu = C$. Thus, $p \notin C$ means $p \notin C\mu$. $\therefore C\mu \subset (C \cup \{p\})\mu$.

($\Leftarrow$). Suppose $p \in U \setminus C$ implies $C\mu \subset (C \cup \{p\})\mu$. We must show that $C = C\mu$.

If $p \in U \setminus C$ then either $p \in C\mu \setminus C$ or $p \in U \setminus C\mu$. Suppose there is a $p \in C\mu \setminus C$. Then by monotonicity and idempotence, $C \cup \{p\} \subseteq C\mu$ and $(C \cup \{p\})\mu \subseteq C\mu$. By our hypothesis this $p$ cannot be. Thus $C\mu \setminus C = \emptyset$, or equivalently $C = C\mu$. □

Note (a) The $\{p\}$ described in Theorem **[CO 7]** is in $Ky\ \mu$. This follows because $C \cup \{p\}$ must contain a key $K$ with $K \setminus \{p\} \subseteq C$ and $p \in K$. That means closed sets can only grow by adding key / independent elements. In matroids and simplicial complexes the independent sets grow when we add $\{p\}$ to an independent set. That is if $C$ is a closed set in a matroid or simplicial complex, and $I$ is an independent subset of $C$ then $I \cup \{p\}$ is independent for every $p \in U \setminus C$. That is not the case with general closures, $I \cup \{p\}$ may be dependent. However, if $I$ is a key of $C$, then $I \cup$

$\{p\}$ will contain a key of $(C \cup \{p\})\mu$ that includes $\{p\}$. Thus, all closed sets in every closure are closed in a very similar way. It takes new independent information to grow.

Note (b) Theorem [**CO 7**] gives the method to construct flats for any hereditary collection $H$. Define $\delta_H: H \to 2^U$ given by $K \mapsto K\delta_H = K \cup \{p \in U \setminus K | K \cup \{p\} \notin H\}$ [R2, 4.2]. The induced closure $\delta_H^+$ using Extend by Closure [**EC 1**] in the next section will produce this closure. A lattice construct called the determined join [R2, 3.2] will produce the same closure function. The closed sets of this extended function are the flats of a matroid or simplicial complex. This method gives only one closure for one hereditary collection. Call this the *flat closure*. There may be many closures associated with one hereditary collection. In the next section [**EC**] we give a bottom-up method to compute this closure and any other closures. In the section [**FL**] we show that the extension by closure and determined join derived from $\delta_H$ are equal.

Before leaving this section, we need some bookkeeping. What is $\emptyset\mu$?

[**CO 8**] *Theorem*: $\emptyset\mu = \{p \in U | \{p\} \notin Ky\ \mu\}$

Proof:

If $p \in \emptyset\mu$ then $\{p\}\mu \subseteq \emptyset\mu$, but $\emptyset \subseteq \{p\}$ means $\emptyset\mu \subseteq \{p\}\mu$. Thus $\{p\}\mu = \emptyset\mu$ and $\{p\}$ cannot be a key of $\mu$.

If $\{p\}$ is not a key of $\mu$ then the only key for $\{p\}\mu$ is $\emptyset$. Thus $p \in \emptyset\mu$. □

## [EC] Extension by Closure -
### Covering Large Spaces Using Small Subsets

In this section we depart from the pattern used so far. Until now we start with a closure $\mu$ and study induced structures like keys and closed sets on its domain. Now, we start with functions that map a finite power set $2^U$ into itself and give a method to extend these functions to closure operators. That is, we give a method to create a finite closure $\mu$. Every finite closure can be derived in this way.

In linear algebra we use a finite set of vectors to compute any member of a possibly infinite span in finite time. This is called extension by linearity.

Here we explore set functions to compute closures. With vectors there are linear combinations. With closure we have three properties, inclusion, monotonicity and idempotence to extend functions on sets to closures. These identify the sets and set operations that can be used to build closed sets from a function.

Let $\alpha$ be a function from $2^U$ into $2^U$. The recurrence relation below uses set operations that are consistent with closure properties to produce the closure operator $\alpha^+: 2^U \to 2^U$ for any such $\alpha$. We call this $\alpha^+$ the extension of $\alpha$ by closure. Just as we never construct an entire vector space, we do not create the entire closure space. Instead, given $\alpha$ and given any $X \subseteq U$, we can derive $X\alpha^+$ in finite time. The recurrence relation below is not optimum. Maier [R1, Algorithm 4.4] gives a linear time algorithm to compute $X\alpha^+$.

The suboptimal recurrence below is useful, however, because it produces an increasing sequence of sets that culminate in the closed set. The recurrence continues but it never changes beyond that point. The increasing sequence of sets are used in this discussion instead of Maier's Derivation DAG's [R1, Definition 4.5].

**[EC 1] The Extension by Closure Recurrence**

1) $X\alpha_0 = X$.
2) For $t \geq 0$, $X\alpha_{t+1} = X\alpha_t \cup \Delta(X\alpha_{t+1})$.
    Where $\Delta(X\alpha_{t+1}) = \cup\{T \in 2^U | S \subseteq X\alpha_t \ \& \ (S,T) \in \alpha\}$.

3) If $X\alpha_{t+1} = X\alpha_t$ then $X\alpha^+ = X\alpha_t$.

Note (a): There is a first $N_X$ such that $X\alpha_n = X\alpha_{N_X}$ for $n \geq N_X$. $N_X \leq |U \setminus X|$. $X\alpha^+ = X\alpha_{N_X}$.

$\alpha^+$ is a closure operator on $2^U$. Before we prove this fact note that the construction is consistent with the closure axioms and the function $\alpha$. For any $S \subseteq U$, $S \subseteq S\alpha^+$ by construction. If we want to extent $\alpha$ to a closure, then if $(S,T) \in \alpha$, the set $T$ should be in $S\alpha^+$ to be consistent with monotonicity. This is true by the construction of the $\Delta$ terms. Finally, we would like idempotence. The terminal value $N_X$ described above will provide idempotence.

[EC 2] *Theorem:* The functions $\alpha^+$ given by $X \mapsto X\alpha^+$ is a closure operator on $2^U$.

Proof:

Inclusion: $X \subseteq X\alpha^+$ by construction.

Monotonicity: Let $X \subseteq Y$. Show that $X\alpha^+ \subseteq Y\alpha^+$.

Prove $X\alpha_n \subseteq Y\alpha_n$ for $n \geq 0$ by induction.

$(n = 0)$ $X\alpha_0 = X \subseteq Y = Y\alpha_0$.

Inductive Hypothesis: For $n > 0$ assume $X\alpha_n \subseteq Y\alpha_n$.

$X\alpha_{n+1} = X\alpha_n \cup \Delta(X\alpha_{n+1})$ and $Y\alpha_{n+1} = Y\alpha_n \cup \Delta(Y\alpha_{n+1})$.

For each $(S,T) \in \alpha$, if $S \subseteq X\alpha_n$ then $S \subseteq Y\alpha_n$ by the induction hypothesis. Thus, if $T$ is used to build $\Delta(X\alpha_{n+1})$ then $T$ is used to build $\Delta(Y\alpha_{n+1})$. $\therefore$ $\Delta(X\alpha_{n+1}) \subseteq \Delta(Y\alpha_{n+1})$. This means: $X\alpha_{n+1} = X\alpha_n \cup \Delta(X\alpha_{n+1}) \subseteq Y\alpha_n \cup \Delta(Y\alpha_{n+1}) = Y\alpha_{n+1}$.

Thus $X\alpha_n \subseteq Y\alpha_n$ for $n \geq 0$ by induction. Specifically, if $n > \max(N_X, N_Y)$ then $X\alpha^+ = X\alpha_n \subseteq Y\alpha_n = Y\alpha^+$.

Idempotence:

We compute $X\alpha^{++}$ using extension by closure:

$X\alpha^{++} = (X\alpha^+)\alpha^+ = (X\alpha_{N_X})\alpha^+$.

$(X\alpha_{N_X})\alpha_0 = X\alpha_{N_X}$.

$(X\alpha_{N_X})\alpha_1 = X\alpha_{N_X} \cup \left(\cup\{T \in 2^U | S \subseteq \alpha_{N_X} \& (S,T) \in \alpha\}\right)$

$= X\alpha_{N_X} \cup \Delta(X\alpha_{N_X+1})$

$= X\alpha_{(N_X+1)} = X\alpha^+$.

$X\alpha^{++} = X\alpha^+$ for every $X \in 2^U$. □

The theorem below shows that the keys paired with their closures are sufficient to determine the entire closure.

**[EC 3] Theorem:** Let $\mu$ be a closure operator on $2^U$. Then the function $\mu|_{Ky\ \mu}$ satisfies $\mu = (\mu|_{Ky\ \mu})^+$.

Proof: We must show that $\mu = (\mu|_{Ky\ \mu})^+$.

Compute $X(\mu|_{Ky\ \mu})^+$ using the algorithm.

$X(\mu|_{Ky\ \mu})_0 = X$, $X(\mu|_{Ky\ \mu})_1 = X \cup \Delta\left(X(\mu|_{Ky\ \mu})_1\right) = X \cup (\cup \{K\mu \mid K \subseteq X \ \& \ (K, K\mu) \in \mu|_{Ky\ \mu}\})$.

Every $K\mu$ in the union above satisfies, $K\mu \subseteq X\mu$ because $K \subseteq X$. Thus, $\Delta\left(X(\mu|_{Ky\ \mu})_1\right) \subseteq X\mu$. However, there is at least one key in the union that is a key of $X$ because we are looking at all keys in $2^X \cap Ky\ \mu$. Thus, $\Delta\left(X(\mu|_{Ky\ \mu})_1\right) = X\mu$, and $X(\mu|_{Ky\ \mu})^+ = X\mu$. □

Note (a) In future we only need to show that $K\alpha^+ = K\mu$ for every $K \in Ky\ \mu$ to show that $\alpha^+ = \mu$.

Every finite closure can be produced from Extension by Closure. Thus, when we refer to a universal closure $\mu$ we need not think of an abstract closure given from thin air. We can think of a $\mu$ constructed from a finite set valued function $\alpha$ that was built to describe the semantic requirements of a model, algebraic or relational. More important, if we have a small $\alpha$ with $\alpha^+ = \mu$, we never need to explicitly enumerate $\mu$. The function $\alpha$ paired with Extend by Closure is $\mu$.

## Functional Dependencies and Closures

When data for an enterprise is stored in tables, certain columns of data determine the values in other columns of data. This results from real constraints of the enterprise being modeled. The fact that {**City**, **Year**} determines {**Rainfall Total**} is one example. Database designers must identify all such dependencies in an enterprise they model.

Formally this is recorded as follows. Let $T$ be a table that includes columns labeled with {**City**, **Year**, **Rainfall Total**}. Let $r_1$ and $r_2$ be any pair of rows in $T$. Each row defines a function from the set of column names into the domain of legitimate values for cells in the table. The fact that {**City**, **Year**} determines {**Rainfall Total**} means that if $r_1|_{\{City,Year\}} = r_2|_{\{City,Year\}}$ then $r_1|_{\{Rainfall\ Total\}} = r_2|_{\{Rainfall\ Total\}}$. The pair ({**City**, **Year**},{**Rainfall Total**}) is a functional dependency on the table $T$.

Database designers collect all the functional dependencies for a table in a binary relation called the family of functional dependencies. Usually these families are not functions, because database designers are only interested in collecting dependency pairs that can be gleaned from careful analysis of the enterprise. For example, ({**City**, **Year**},{**Mean Income** }) is another dependency in the family. Extend by closure also works for families of functional dependencies even if they are not presented as functions.

However, for any binary relation $\rho$ on the subsets in $2^U$, there is one function $\alpha$ with the same domain as $\rho$ such that $\rho^+ = \alpha^+$ and $\alpha \subseteq \alpha^+$. Given a closure $\mu$ we always use the unique $\alpha \subseteq \mu$ with $\alpha^+ = \mu$ to represent all binary relations $\rho$ with the same domain that satisfy $\rho^+ = \mu$. This does restrict the generality that is useful for the database designer's manual data analysis, but it does not limit the set of closures that can be expressed.

## Spans, Covers and Nonredundant Functions

The function $\mu$ denotes a closure operator on $2^U$. The symbols $\alpha$, $\beta$ and $\gamma$ denote subsets of $\mu$. $\alpha^+$, $\beta^+$ and $\gamma^+$ are the closures generated with Extend by Closure.

**[EC 4]** Span: The *span* of $\alpha$ is $\alpha\mu^* = \alpha^+ \cap \mu$. That is, all pairs $(S, S\alpha^+)$ in $\alpha^+$ such that $S\alpha^+ = S\mu$. $\mu^*$ is called the *span operator*.

**[EC 5]** Cover: The function $\alpha \subseteq \mu$ covers $\mu$ when $\alpha^+ = \mu = \alpha\mu^*$. That is, the span of $\alpha$ is all $\mu$.

**[EC 6]** Nonredundant / Redundant: A cover $\alpha$ is *nonredundant* whenever no proper subset of $\alpha$ has the closure $\alpha^+$. A a subset of a nonredundant cover $\alpha$ is called a *nonredundant function*. Any cover or function $\alpha$ is *redundant* if there is a $(S, S\mu) \in \alpha$ such that $(\alpha \setminus \{(S, S\mu)\})^+ = \alpha^+$, or equivalently, $(S, S\mu) \in (\alpha \setminus \{(S, S\mu)\})^+$.

**[EC 7]** *Lemma*: If a cover $\alpha$ includes a redundant subset $\beta$, then $\alpha$ is redundant.

Proof: Let $\alpha$ be a nonredundant cover. $\alpha = \beta \cup \gamma$ where $\gamma = \alpha \setminus \beta$. There is an $\eta = \beta \setminus \{(X, X\mu)\}$ for some $(X, X\mu) \in \beta$ such that $\eta^+ = \beta^+$ because $\beta$ is redundant. $\alpha \setminus \{X, X\mu\} = \eta \cup \gamma$. Now, $\alpha^+ \subseteq (\beta^+ \cup \gamma^+)^+ = (\eta^+ \cup \gamma^+)^+ \subseteq (\eta \cup \gamma)^+ = (\alpha \setminus \{(X, X\mu)\})^+ \subseteq \alpha^+$. Hence $\alpha$ is redundant. □

To produce a nonredundant cover from a redundant cover, just remove redundant pairs one by one terminating when no redundant pair can be found. This can be done in one pass over the function, [R1, Algorithm 5.4].

**[EC 8]** *Theorem*: The span operator $\mu^*: 2^\mu \to 2^\mu$ given by $\alpha \mapsto \alpha\mu^*$ is a closure operator on $2^\mu$, and the keys of $\mu^*$ are the nonredundant functions $\alpha \in 2^\mu$.

Proof: (Closure)

*Inclusion*: Pick $\alpha$ in $2^\mu$. Then $\alpha \subseteq \mu$. Hence, $\alpha \subseteq \alpha^+$ because $X\mu = X\alpha \subseteq X\alpha^+ \subseteq X\mu$ for every $(X, X\mu) \in \alpha$. Thus $\alpha \subseteq \alpha^+ \cap \mu = \alpha\mu^*$.

*Monotonicity*: If $\alpha \subseteq \beta$ then $X\alpha^+ \subseteq X\beta^+$ for any $X \subseteq U$. This is because any set $S\mu$ used from a pair in $\alpha$ to produce $X\alpha^+$ will also be used to produce $X\beta^+$. Every $S\mu$ used in Extend by Closure is a subset of $X\mu$. That means $X\alpha^+ = X\mu \Rightarrow X\beta^+ = X\mu$. $\therefore \alpha\mu^* = \alpha^+ \cap \mu \subseteq \beta^+ \cap \mu = \beta\mu^*$.

*Idempotence*: We must show that $\alpha\mu^*\mu^* = \alpha\mu^*$. We have shown that: $\alpha \subseteq \alpha\mu^*$. Also, $\alpha\mu^* \subseteq \alpha^+$ by definition. Applying the monotonicity of $\mu^*$ (just proved) gives $\alpha\mu^* \subseteq \alpha\mu^*\mu^* \subseteq \alpha^+\mu^*$. Now, $\alpha^+\mu^* = \alpha^{++} \cap \mu = \alpha^+ \cap \mu = \alpha\mu^*$. Thus $\alpha\mu^*\mu^* = \alpha\mu^*$.

($Ky\ \mu^*$ equals the nonredundant functions $\alpha \subseteq \mu$)

Let $\alpha \in Ky\ \mu^*$ be a maximal subset regarding inclusion. No proper subset of $\alpha$ has the same span $\alpha\mu^* = \mu$. Let $\beta = \alpha \setminus \{(S, S\mu)\}$ for an arbitrary $(S, S\mu) \in \alpha$. $\beta\mu^* \subset \alpha\mu^*$. Expanding, $\beta^+ \cap \mu \subset \alpha^+ \cap \mu$. This says there must be some $(T, T\mu) \in \alpha^+$ that is not in $\beta^+$. Thus $\alpha$ is nonredundant because no proper subset has the same closure. Since $Ky\ \mu^*$ is hereditary, any subset $\beta$ of $\alpha$ is in $Ky\ \mu^*$, and $\beta \subseteq \alpha$ is nonredundant, by **[EC 7]**. This is true for any maximal $\alpha \in Ky\,\mu^*$ and $\beta \subseteq \alpha$. Thus $Ky\,\mu^*$ contains all nonredundant functions.

$\alpha \subseteq \mu$ nonredundant means that $\alpha$ is a subset of a nonredundant cover $\beta$. It follows that $(S, S\mu) \notin (\alpha \setminus \{(S, S\mu)\})^+$ for any $(S, S\mu) \in \alpha$. Otherwise, $\alpha$ and $\beta$ are redundant, **[EC 7]**. $\therefore (\alpha \setminus \{(S, S\mu)\})^+ \cap \mu \neq \alpha^+ \cap \mu$ because $(S, S\mu)$ is not in the left-hand side. This is true for every $(S, S\mu) \in \alpha$. That means no proper subset of $\alpha$ has the same span $\alpha\mu^*$. $\alpha \in Ky\,\mu^*$. □

The span operator $\mu^*$ will prove useful when comparing covers in the next sections. It is also a flat-closure as we prove next.

## [FL] The Flat Closure

The canonical closure associated with every hereditary collection is the flat closure. This is the closure for matroids and simplicial complexes.

**[FL 1]** Dependency Function $\delta_H$: The *dependence function* of a hereditary collection $H$ is the function $\delta_H: H \to 2^U$ given by $K\delta_H = K \cup \{p \in U \setminus K \mid \{p\} \cup K \notin H\}$.

**[FL 2]** The *flat closure* of $H$, denoted $\kappa_H: 2^U \to 2^U$, is defined as follows: $X\kappa_H = \cap \{I\delta_H \mid I \in H\ \&\ X \subseteq I\delta_H\}$. [R2, 4.2]

To verify that this definition make sense, we show $\kappa_H$ is a closure and that that $\delta_H^+ = \kappa_H$. That is the bottom-up extension by closure of $\delta_H$ equals the top-down definition of the flat closure $\kappa_H$. The top-down definition is the determined join of lattice theory [R2, 3.2]. This is the first time in this exposition a closure operator has been defined from its set of keys. Up until now we have created hereditary collections (keys) from closures that were generated bottom up using Extend by Closure.

**[FL 3]** Ancestors: For any $X \in 2^U$ the set of *ancestors of $X$ regarding $\delta_H$* is given by: $An(X, H) = \{I\delta_H \mid I \in H\ \&\ X \subseteq I\delta_H\}$.

Note (a) $I\delta_H = I\kappa_H = \cap\, An(I\delta_H, H)$ for $I \in H$.

**[FL 4]** *Lemma*: The function $\kappa_H$ is a closure with key set $H$.

Proof:

(Inclusion): This follows from the definition.

(Monotonicity): $X \subseteq Y \Rightarrow An(Y, H) \subseteq An(X, H)$. That means $X\kappa_H \subseteq Y\kappa_H$ because $\cap\, An(X, H) \subseteq \cap\, An(Y, H)$.

(Idempotence: $X\kappa_H = X\kappa_H\kappa_H$):

Claim (a) $An(X, H) = An(X\kappa_H, H)$:

$X \subseteq X\kappa_H \Rightarrow An(X\kappa_H, H) \subseteq An(X, H)$.

To see that $An(X, H) \subseteq An(X\kappa_H, H)$, let $I\delta_H \in An(X, H)$. Then $X \subseteq I\delta_H$ implies $X\kappa_H \subseteq I\delta_H$ because $X\kappa_H$ is the intersection of all such $I\delta_H$. Thus $I\delta_H \in An(X\kappa_H, H)$. Claim (a) is proved.

$X\kappa_H = \cap\, An(X, H) = \cap\, An(X\kappa_H, H) = X\kappa_H\kappa_H$.

Thus $\kappa_H$ is a closure on $2^U$.

(Key Set: $Ky\,\kappa_H = H$):

Suppose $I \in H$. $H$ is hereditary so $I \setminus \{p\} \in H$ and $(I \setminus \{p\})\delta_H = (I \setminus \{p\})\kappa_H$. That set does not include $p$ because $I \in H$. Thus, no proper subset of $I$ has the same closure $I\kappa_H$. Hence, $I$ is a key of $\kappa_H$. $I \in Ky\,\kappa_H$.

Suppose $K \in Ky\,\kappa_H$. Let $I$ be a maximal independent set in $K$. For every $p \in K \setminus I$ the set $I \cup \{p\} \notin H$. Thus $K \subseteq I\delta_H$. This says $I\delta_H \in A(X, H)$ and $K\kappa_H = \cap\, A(X, H) \subseteq I\delta_H = I\kappa_H$. That is $K\kappa_H \subseteq I\kappa_H$. However, $I \subseteq K$ means that $I\kappa_H \subseteq K\kappa_H$. We have $I\kappa_H = K\kappa_H$. If $I$ is a proper subset of $K$ then $K$ is not a key of $K\kappa_H$, contradicting our choice of $K$. Hence, $I = K$, and $K \in H$. □

**[FL 5]** *Corollary*: If $H$ is a hereditary collection over $U$ and $X \in 2^U$ then $X\kappa_H = I\delta_H$ for every maximal independent set $I \subseteq X$.

Proof: Follows from $Ky\,\kappa_H = H$ and $I\delta_H = I\kappa_H$. $I\kappa_H \subseteq X\kappa_H$ by monotonicity, and $X\kappa_H \subseteq I\delta_H = I\kappa_H$ because $I$ is maximal independent in $X$. □

**[FL 6]** *Theorem*: $\delta_H^+ = \kappa_H$.

Proof: $\delta_H = \{(I, I\delta_H) | I \in H\}$. $H = Ky\,\kappa_H$ and $I\delta_H = I\kappa_H$ for $I \in H$. Equivalently, $\delta_H = \kappa_H|_{Ky\,\kappa_H}$. Theorem **[EC 3]** says $(\mu|_{Ky\,\mu})^+ = \mu$ for any closure $\mu$. Thus, $\delta_H^+ = \kappa_H$ □

The top-down and bottom-up constructions of the flat closure agree. The following corollary shows that the span operator, $\mu^*$, is a flat closure.

**[FL 7]** *Corollary*: Let $H$ be the set of all nonredundant functions $\alpha \in 2^\mu$, where $\mu$ is closure operator on $2^U$. The span operator $\mu^*: 2^\mu \to 2^\mu$ is the flat closure for $H$.

Proof: Theorem **[EC 8]** showed that $\mu^*$ is a closure with key set $H$, the set of all nonredundant functions. Let $\alpha \in H$. Let $p = (X, X\mu)$

If $p \in \alpha\mu^* \setminus \alpha$ then $\alpha' = \alpha \cup \{p\}$ is redundant. That is because $p \in (\alpha' \setminus \{p\})^+ = \alpha^+$. Thus, $\alpha\mu^* \subseteq \alpha\delta_H$.

If $p = (X, X\mu) \in \alpha\delta_H$ then $\alpha \cup \{p\} \notin H$. The hereditary collection $H$ is the set of all nonredundant functions. Thus, $\alpha \cup \{p\}$ is redundant. By definition, $(\alpha \cup \{p\})^+ = \alpha^+$. Thus, $p = (X, X\mu) \in \alpha^+ \cap \mu = \alpha\mu^*$. ∴ $\alpha\delta_H \subseteq \alpha\mu^*$.

We have $\mu^*|_H = \delta_H$ which means $\mu^* = \delta_H^+ = \kappa_H$, the flat closure for $H$, **[FL 6]**. □

# *[MAT] The Matroid $\mathcal{M}(\mu)$*

## Maier's Bijection

We demonstrated that the span operator $\mu^*$ of any finite closure $\mu$ is a closure for the hereditary collection $H$ composed of the nonredundant functions $\alpha \in 2^\mu$, [**EC 8**]. Moreover, $\mu^*$ is the flat closure for this hereditary collection [**FL 7**].

Thus, every finite closure $\mu$ induces a simplicial complex with ground set $\mu$ and closure $\mu^*$. In this section we show that this simplicial complex is a matroid denoted by $\mathcal{M}(\mu)$ and that the nonredundant covers $\alpha$ of $\mu$ form the bases of that matroid.

We will do this by demonstrating a bijection between nonredundant covers that also supports exchange. This construction was introduced by Maier [**R1**]. Maier's goal was to devise an algorithm that would find minimal covers for families of functional dependencies that could be used to partition a database into a set of tables that eliminated avoidable data repetition and preserved the family of functional dependencies for the model. The algorithm created an efficient representation of the enterprise that was faithful to the database designer's analysis. As a side effect, Maier showed that his set of minimal covers defined the bases of a matroid although Maier did not claim this in [**R1**].

Maier used families of functional dependencies that were not necessarily functions. Maier needed an extra step to compact his families into minimal functions. We will only consider functions $\alpha \subseteq \alpha^+$. This does not restrict the closure operators that are covered using Extend by Closure, but it does ensure that $\alpha$ nonredundant means $\alpha$ is minimal.

General families of dependencies are more reasonable for database designers because there may be thousands of such dependencies to consider in a large enterprise. It is best to use an algorithm to put them into a regular form. Maier describes this process completely [R1, Ch. 1–5].

The more restricted form $\alpha \subseteq \alpha^+$ is used here because it makes proving the matroid property more direct. We will discuss Maier's bijection, direct determination, below, but first we need some notation and one critical structure theorem, [**MAT 4**].

## Range Restrictions

*Theorem* [CO 7] demonstrates an onion structure for closed sets. Given a closed set $C$ we grow it only by layering new key elements from outside the set to produce a new closed set that is proper superset. The Extension by Closure recurrence uses this structure to compute closed sets from below. The set $X\alpha^+$ only requires pairs of the form $(S, S\mu) \in \alpha$ that satisfy $S, S\mu \subseteq X\alpha^+$, [EC 1]. Each $X\alpha_n$ in the recurrence grows the set until a closed set boundary is reached. Nothing from inside can grow the set further. That explains the idempotence property. The sets described below will assist in discussing the range restrictions that bound this process.

[**MAT 1**] Let $C$ be a closed set of $\mu$ and $\alpha \subseteq \mu$. Consider the following range restrictions of $\alpha$ regarding closed set $C$.

  Body:     $\alpha_{\subseteq C} = \{(S, S\mu) \in \alpha |\ S\mu \subseteq C\}$
  Interior: $\alpha_{\subset C} = \{(S, S\mu) \in \alpha |\ S\mu \subset C\}$.
  Top:      $\alpha_{=C} = \{(S, S\mu) \in \alpha | S\mu = C\}$.

These restrictions are called the *body*, *interior* and *top* of $\alpha$ regarding $C$ respectively. The body is all $\alpha|_{2^C}$. The interior is the (set, closure) pairs that do not reach $C$. The top is the boundary between the interior and exterior pairs $(S, S\mu)$ where $S\mu \subseteq U\setminus C$.

This non-standard notation is useful because it describes limits to the range of functions rather than the domain. The notation for domain restriction is clumsy for restricting the range.

Note (a) $(\alpha_{\subseteq C})^+$ and $(\alpha_{\subset C})^+$ are denoted by $\alpha^+_{\subseteq C}$ and $\alpha^+_{\subset C}$ respectively. $\alpha^{++}_{\subset C} = ((\alpha_{\subset C})^+)^+$, etc. In general, $\alpha^+_{\subset C} \neq (\alpha^+)_{\subset C}$. The closure of the interior can include elements of the top, but the interior of the closure cannot. Equality only holds when the top is empty.

Note (b): The closure universe used for many functions will be the body $\mu_{\subseteq C}$ rather than the full $\mu$. Again, this works like an onion. $\mu_{\subseteq C}$ is a closure operator that agrees with $\mu$ on $2^C$. Denote: $\mu^*_{\subset C} = (\mu_{\subset C})^*$ and $\mu^*_{\subseteq C} = (\mu_{\subseteq C})^*$.

**[MAT 2]** *Lemma*: If $\alpha \subseteq \mu$ covers $\mu$ then

(1) $\alpha^+_{\subseteq C} = \mu_{\subseteq C}$ and

(2) $\alpha^+_{\subset} \subseteq \mu_{\subseteq C}$.

Proof: (1) The domain of $\mu_{\subseteq C}$ is $2^C$, and for every $X \subseteq C$, $X\mu = X\alpha^+ = X\alpha^+_{\subseteq C}$ because every pair used to compute $X\alpha^+$ resides in $\alpha_{\subseteq C}$ [EC 1].

(2) $X\alpha^+_{\subset} = X\alpha^+_{\subseteq C}$ for every $X$ in the domain of $\alpha_{\subset C}$. Thus, $\alpha^+_{\subset C} \subseteq \alpha^+_{\subseteq C} = \mu_{\subseteq C}$. □

**[MAT 3]** *Lemma*: Let $\alpha \subseteq \mu$ and $\beta \subseteq \mu$ cover $\mu$, a closure operator on $2^U$. For any closed set $C$, if $\alpha_{=C} = \emptyset$ then $\beta^+_{\subset C} = \beta^+_{\subseteq C}$.

Proof: $\alpha_{\subset C} \subseteq \beta^+_{\subset C}$, because $X\alpha^+ = X\beta^+ = X\mu$ for every $X \subset C$ and $\alpha_{\subset C} \subseteq \alpha^+_{\subset C}$.

We now apply the span operator $\mu^*_{\subseteq C}$. $\alpha_{\subset C}\mu^*_{\subseteq C} \subseteq \beta^+_{\subset C}\mu^*_{\subseteq C}$. Expanding gives $\alpha_{\subset C}\mu^*_{\subseteq C} = \alpha^+_{\subset C} \cap \mu_{\subseteq C} = \alpha^+_{\subset C}$ and $\beta^+_{\subset C}\mu^*_{\subseteq C} = \beta^{++}_{\subset C} \cap \mu_{\subseteq C} = \beta^+_{\subset C}$ by [**MAT 2**](2) and $\beta^{++}_{\subset C} = \beta^+_{\subset C}$, [EC 2]. Hence $\alpha^+_{\subset C} \subseteq \beta^+_{\subset C}$.

But $\alpha_{=C} = \emptyset$ which means $\alpha_{\subset C} = \alpha_{\subseteq C}$. Thus, $\alpha^+_{\subset C} = \alpha^+_{\subseteq C} = \mu_{\subseteq C}$ by [**MAT 2**](1). Concluding, $\mu_{\subseteq C} = \alpha^+_{\subset C} \subseteq \beta^+_{\subset C} \subseteq \beta^+_{\subseteq C} = \mu_{\subseteq C}$, and $\therefore \beta^+_{\subset C} = \beta^+_{\subseteq C}$. □

This is an important lemma because it says that if you have one cover that is empty when restricted to a closed set, then every cover that is nonempty on that closed set is redundant.

[**MAT 4**] *Theorem*: If $\alpha \subseteq \mu$ and $\beta \subseteq \mu$ are nonredundant covers of $\mu$ a closure on $2^U$ and $C$ is a closed set of $\mu$. Then $\alpha_{=C} = \emptyset$ if and only if $\beta_{=C} = \emptyset$.

Proof: Without loss of generality, assume the top set $\alpha_{=C} = \emptyset$. From the *Lemma* [**MAT 3**], $\beta^+_{\subset C} = \beta^+_{\subseteq C}$. That means $(\beta_{\subseteq C} \setminus \beta_{=C})^+ = \beta^+_{\subseteq C}$. If $\beta_{=C} \neq \emptyset$ then $\beta$ is redundant. □

Every closed set $C$ of a closure $\mu$ can be marked as either redundant or nonredundant depending on whether the corresponding top $\alpha_{=C}$ is empty for just one nonredundant cover. There are many such closed sets. For example, every closed set $C$ that is also independent (a key) has an empty top set in every nonredundant cover. This is because the only member of the top has the form $(C, \mu C) = (C, C)$. Reflexive pairs are redundant in every cover that contains them.

## Direct Determination

The interior set $\alpha_{\subset C}$ can determine the top set $\alpha_{=C}$ as witnessed by *Lemma* [**MAT 3**]. This section will explore determination from the interior of a cover more carefully. It will be useful for identifying redundancy, linking members of different nonredundant covers and facilitating exchange. Maier [R1] identified this relation and called it *direct determination*. Most of the discussion that follows was developed by Maier. It is formulated quite differently, but parallel concepts are cited. Maier's symbol → denotes the relation.

**[MAT 5]** Direct Determination: The set $X$ determines $Y$ directly, denoted by $X \dashrightarrow Y$, whenever $X\mu = Y\mu = C$ and $Y \subseteq X\mu^+_{\subset C}$. [R1, Definition 5.9]

Direct determination says we can determine that $Y \subseteq X\mu$ without consulting any top element $(Z, C) \in \mu$. The restriction $\mu_{\subset C}$ and the one top element $(Y, C)$ are sufficient to derive $X\mu = C$ using extension by closure.

**[MAT 6]** Direct determination is reflexive, projective, and transitive.

Reflexivity: $X \subseteq X\mu^+_{\subset C}$ by the definition of $X\mu^+_{\subset C}$ [EC 1]. Thus, $X \dashrightarrow X$.

Projectivity: If $Y \subseteq X$ and $X\mu = Y\mu = C$ then $X \dashrightarrow Y$ because $Y \subseteq X \subseteq \mu^+_{\subset C}(X)$.

Transitivity: Suppose $X \dashrightarrow Y$ and $Y \dashrightarrow Z$. Then $Y \subseteq X\mu^+_{\subset C}$ and $Z \subseteq Y\mu^+_{\subset C}$. $Y \subseteq X\mu^+_{\subset C}$, gives $Y\mu^+_{\subset C} \subseteq X\mu^{++}_{\subset C} = X\mu^+_{\subset C}$. Hence $Z \subseteq Y\mu^+_{\subset C} \subseteq X\mu^+_{\subset C}$. □

Direct determination is also symmetric for members of nonredundant covers. In this special case, direct determination defines an equivalence relation, with class members that are interchangeable elements of nonredundant covers.

## Nonredundant Covers as Bases

We begin by sharpening the result in *Lemma* **[MAT 2]**(2). The bottom-up proof of this lemma is characteristic of direct determination proofs.

**[MAT 7]** *Lemma*: If $\alpha \subseteq \mu$ is a cover of $\mu$, a closure over $2^U$, then $\alpha^+_{\subset C} = \mu^+_{\subset C}$. ([R1] *Lemma* 5.5)

Proof: ($\alpha^+_{\subset C} \subseteq \mu^+_{\subset C}$): $\alpha_{\subset C} \subseteq \mu_{\subset C}$, because the $\alpha \subseteq \mu$. Applying the closure $\mu^*_{\subseteq C}$ gives $\alpha_{\subset C}\mu^*_{\subseteq C} \subseteq \mu_{\subset C}\mu^*_{\subseteq C}$. Expanding, $\alpha^+_{\subset C} \cap \mu_{\subseteq C} \subseteq \mu^+_{\subset C} \cap \mu_{\subseteq C}$. Both are subsets of $\mu_{\subseteq C}$ [MAT 2](2), so $\alpha^+_{\subset C} \subseteq \mu^+_{\subset C}$.

($\mu^+_{\subset C} \subseteq \alpha^+_{\subset C}$): We first show $\mu_{\subset C} \subseteq \alpha^+_{\subset C}$. To see this note: For $X\mu \subset C$, apply Extend by Closure [EC 1] to compute $X\alpha^+ = X\mu$. Let $X =$

$X\alpha_0, \ldots, X\alpha_{N_X} = X\mu$ be the sequence to compute $X\mu$. Each $X\alpha_i \subset C$. Thus, every pair $(S, S\mu) \in \alpha$ used to compute each $\Delta(X\alpha_i)$ must satisfy $S, S\mu \subset C$. That means every $(S, S\mu)$ used is an element of $\alpha_{\subset C}$, or equivalently $X\alpha^+ = X\alpha^+_{\subset C}$ for every pair, $X, X\mu \subset C$. $\therefore \mu_{\subset C} \subseteq \alpha^+_{\subset C}$.

Applying $\mu^*$ gives $\mu_{\subset C}\mu^* \subseteq \alpha^+_{\subset C}\mu^*$ and expanding $\mu^+_{\subset C} \cap \mu_{\subseteq C} \subseteq \alpha^{++}_{\subset C} \cap \mu_{\subseteq C}$. $\alpha^+_{\subset C} = \alpha^{++}_{\subset C}$, and $\mu^+_{\subset C}, \alpha^+_{\subset C} \in \mu_{\subseteq C}$, gives $\mu^+_{\subset C} \subseteq \alpha^+_{\subset C}$. □

The preceding lemma means that if we can show that $Y \subseteq X\alpha^+_{\subset C}$ for any cover $\alpha$ we have $X \rightarrow Y$.

[**MAT 81**] *Lemma*: Let $\alpha \subseteq \mu$ be a cover of the closure $\mu$ over $2^U$. If

(1)      $C$ be a closed set,

(2)      $Y\mu = C$, and

(3)      $\alpha_{=C} \neq \emptyset$,

then there is a $Z \in \alpha^{-1}(C)$ such that $Y \rightarrow Z$. (Maier[R1] *Lemma* 5.7)

Proof:

Case (A): There is a $Z \in \alpha^{-1}(C)$ with $Z \subseteq Y$: Then $Y \rightarrow Z$ by projectivity [**MAT 6**].

Case (B): $Y \rightarrow C$: In this case $Y \rightarrow Z$ for every $Z \in \alpha^{-1}(C)$ by transitivity [**MAT 6**].

Case (C): $Y$ has no subset in $\alpha^{-1}(C)$ and $Y$ does not determine $C$ directly.

This ensures $Y$ is a proper subset of $C$ and every proper subset $S \subset Y$ satisfies $S\mu \subset C$. This case assumes $Y$ does not determine $C$ directly. $\therefore Y$ must use a pair in $\alpha_{=C}$ to compute $Y\alpha^+ = C$, because $Y\alpha^+ = C$.

Applying Extend by Closure [EC 1] to $Y$ yields the following: $Y\alpha_{N_Y-1}$ must contain an element $Z \in \alpha^{-1}(C)$ so that $(Z, C)$ can be used to compute $Y\alpha_{N_Y} = C$. No earlier term of the sequence $\{Y\alpha_k\}_{k=1}^{N_Y}$ can contain a member of $\alpha^{-1}(C)$, or $Y\alpha^+ = C$ would be encountered earlier. Stated differently, $Y\alpha_{N_Y-1} \subseteq Y\alpha^+_{\subset C}$ because every pair used to compute $Y\alpha_{N_Y-1}$ is taken from $\alpha_{\subset C}$. By [MAT 7], $Y \dashrightarrow Y\alpha_{N_Y-1}$ and $Y \dashrightarrow Z$ because $Z \subseteq Y\alpha_{N_Y-1}$. The set $Z$ is the one we seek. □

Direct determination enables exchange between elements of nonredundant covers. It also detects redundancy.

[MAT 9] *Lemma*: Let $\alpha, \beta \subseteq \mu$ be nonredundant covers for the closure $\mu$ over $2^U$. Suppose there is $X \neq Y$ with $(X, C) \in \alpha$, $(Y, C) \in \beta$ such that $X \dashrightarrow Y$. Then

(1) $\left((\alpha \setminus \{(X, C)\}) \cup \{(Y, C)\}\right)^+ = \alpha^+$.

(2) If $\alpha = \beta$ then $(\alpha \setminus \{(X, C)\})^+ = \alpha^+$.

[R1], *Theorem* 5.1, *Lemma* 5.8

Proof:

(1) Let $\gamma = (\alpha \setminus \{(X, C)\}) \cup \{(Y, C)\}$. $\gamma \subseteq \mu$ and $\gamma_{\subset C} = \alpha_{\subset C}$ by construction. $\alpha$ is a cover and $X \dashrightarrow Y$ implies $Y \subseteq X\alpha^+_{\subset} = X\gamma^+_{\subset C}$, [MAT 7]. Consider the derivation of $X\gamma^+$: $X\gamma_0, \ldots X\gamma_{N_X} = X\gamma^+$. The set $X\gamma_{N_X-1}$ will include a $Z \in \gamma^{-1}(C)$. This is guaranteed because $(Y, C) \in \gamma$ and $Y \subseteq X\gamma^+_{\subset C}$. $Z$ may or may not equal $Y$. Every pair used to compute $X\gamma^+_{\subset C}$ is also used to compute $X\gamma^+$ [EC 1]. Thus $(X, C) \in \gamma^+$ and $\alpha \subseteq \gamma^+$. If we apply $\mu^*$ and $\mu = \alpha^+ \cap \mu = \alpha\mu^* \subseteq \gamma^+\mu^* = \gamma^{++} \cap \mu = \gamma^+ \cap \mu \subseteq \mu$. Thus, $\mu = \alpha^+ = \gamma^+$.

(2) If $\alpha = \beta$ then $(\alpha \setminus \{(X, C)\}) \cup \{(Y, C)\} = (\alpha \setminus \{(X, C)\})$. This is because $X \neq Y$ and $(Y, C) \in \alpha = \beta$. Thus, using (1), $(\alpha \setminus \{(X, C)\})^+ = \alpha^+ = \mu$. □

Note (a) The roles of $\alpha$ and $\beta$ are interchangeable. That is $(\beta \setminus \{(Y, C)\}) \cup \{(X, C)\}$ is also a cover.

Note (b): Restating (2) says: If $\alpha$ is a cover and there are distinct $(X, C)$ and $(Y, C)$ in $\alpha$ with $X \dashrightarrow Y$ then $\alpha$ is redundant.

**[MAT 10]** *Lemma*: If $\alpha$ and $\beta$ are nonredundant covers of $\mu$, a closure operator over $2^U$, then direct determination defines a self-inverse bijection between the elements of $\alpha$ and $\beta$.

Proof: Let $C$ be a closed set of $\mu$. If $(X, C) \in \alpha$ there is a $(Y, C) \in \beta$ such that $X \dashrightarrow Y$. This follows from **[MAT 4]** ($\beta_{=C} \neq \emptyset$), and $\beta$ is a cover with $X\mu = X\alpha = C$ (given). The three conditions of **[MAT 8]** are met so such an $(Y, C) \in \beta$ exists. Similarly, there is a $(Z, C) \in \alpha$ such that $Y \dashrightarrow Z$. Transitivity of direct determination gives $X \dashrightarrow Z$ and the nonredundancy of $\alpha$ dictates $X = Z$, *Lemma* **[MAT 9**(2)**]**. We have shown that $\dashrightarrow$ is a self-inverting relation from $\alpha$ to $\beta$.

If $X \dashrightarrow W$ for some $(W, C) \in \beta$ then $W = Y$ by the same reasoning: $Y \dashrightarrow X$, $X \dashrightarrow W \Rightarrow Y \dashrightarrow W$. $\beta$ is nonredundant means $W = Y$. Thus, $\dashrightarrow$ is a function form $\alpha$ into $\beta$.

∴ If $C$ is a closed set of $\mu$, and $\alpha$ and $\beta$ are nonredundant covers of $\mu$ then direct determination is an injection that is its own inverse. The same argument holds for mapping elements of $\beta$ into alpha.

Direct determination is a bijection between nonredundant functions in $2^\mu$. □

**[MAT 11]** *Theorem*: If $\mu$ is a closure over $2^U$ then $\mathcal{M}(\mu) = (\mu, Ky\,\mu^*)$ is a matroid with ground set $\mu$ and independent sets $Ky\,\mu^*$, the non-redundant functions in $2^\mu$.

Proof. We show that the nonredundant covers of $\mu^*$ are the basis sets for $\mathcal{M}(\mu)$ by demonstrating equal cardinality and exchange.

(Equal Cardinality) The bijection, direct determination, between the nonredundant covers just proved in *Lemma* [**MAT 10**], gives equal cardinality.

(Exchange) Let $\alpha \neq \beta$ be nonredundant covers of $\mu$. There is a closed set $C$ such that $\alpha_{=C} \neq \beta_{=C}$. That means both restrictions are nonempty by [**MAT 4**].

Let $(X, C) \in \alpha \setminus \beta$. There is a $(Y, C) \in \beta$ such that $Y \to X$. The pair $(Y, C)$ cannot be in $\alpha$ or $\alpha$ would be redundant. By *Lemma* [**MAT 9**(1)], $\gamma = (\alpha \setminus \{(X, C)\}) \cup \{(Y, C)\}$ is a cover. And $|\gamma| = |\alpha|$ means $\gamma$ is nonredundant because no proper subset can cover $\mu$.

The nonredundant covers are basis sets of the matroid $\mathcal{M}(\mu) = (\mu, Ky\, \mu^*)$. □

The independent sets of this matroid can be quite sparce. The following result gives a description of the singleton sets that are in no independent set.

We denote singleton functions of $\mu^*$ as follows: for each $(X, C) \in \mu$, the singleton function $\pi_{X,C} = \{(X, C)\}$. This can get tricky. $\pi_{X,C}$ is a function. Thus if $X \subseteq U$ the ground set of $\mu$ then $X\pi_{X,C}\mu^*$ is a set in $U$ and $\pi_{X,C}\mu^*$ is the function $\pi_{X,C}^+ \cap \mu$.

Other subsets of $\mu$ are denoted with lower case Greek letters as before.

[**MAT 12**] *Lemma*: If $X \to X\mu$ or $X \notin Ky\, \mu$ then $\pi_{X,X\mu} \in \emptyset\mu^*$ (equivalently $\pi_{X,X\mu} \notin Ky\, \mu^*$).

Proof: Throughout this proof: Let: $C = X\mu = K\mu$ and $\alpha$ be an arbitrary basis, where $K \subseteq X \subseteq C$ and $K$ is a key of $C$ in $\mu$.

(Case 1) Suppose $X \to C$. Transitivity implies $X \to Y$ for every $(Y, C) \in \alpha_{=C}$, [MAT 6].

If $\alpha_{=C} \ne \emptyset$ then $\alpha$ would be redundant if it included $\pi_{X,C}$, [MAT 9](2). Thus $\pi_{X,C} \not\subseteq \alpha$. If $\alpha_{=C} = \emptyset$ then $\pi_{X,C} \not\subseteq \alpha$ because $X\mu = C$.

Thus, if $X \dashrightarrow C$ then $\pi_{X,C}$ cannot be a subset of any nonredundant cover of $\mu$. Thus $\pi_{X,C} \notin Ky\ \mu^*$. Equivalently $\pi_{X,C} \in \emptyset\mu^*$, by [CO 8].

(Case 2) Suppose $(X \notin Ky\ \mu)$, then $K \subset X$ and $\eta = \{(X,C),(K,C)\}$ is redundant. This is true because $X \dashrightarrow K$, [MAT 6]. $\pi_{X,C} \in (\eta \setminus \pi_{X,C})^+ = \pi_{K,C}^+$, [MAT 9](2). There are two cases:

(Case A) $\pi_{K,C} \in \emptyset\mu^*$: This says $\pi_{K,C}\mu^* \in \emptyset\mu^*$. $\pi_{X,C} \subseteq \pi_{K,C}^+ \cap \mu = \pi_{K,C}\mu^*$. $\therefore$ $\pi_{X,C} \in \emptyset\mu^*$.

(Case B) $\pi_{K,C} \in Ky\ \mu^*$: First note, $\pi_{K,C} \notin \pi_{X,C}\mu^*$. This is because $K\pi_{X,C}^+ = K \ne C$ using Extend by Closure and $(K,K) \notin \mu$. This means that $\pi_{X,C}^+ \cap \mu = \pi_{X,C}^+\mu^*$ cannot include $\pi_{K,C}$.

If $\pi_{X,C} \in Ky\ \mu^*$ then $\pi_{X,C}\mu^* = \pi_{X,C}\delta_{Ky\ \mu^*}$ by [**FL 5**]. This says that $\pi_{X,C} \cup \pi_{K,C} \in Ky\ \mu^*$ because $\pi_{K,C} \notin \pi_{X,C}\mu^* = \{\pi_{S,S\mu}|\pi_{X,C} \cup \pi_{S,S\mu} \notin Ky\ \mu^*\}$.

But we already showed that $\eta = \pi_{X,C} \cup \pi_{K,C}$ is redundant. Thus, $\pi_{X,C} \notin Ky\ \mu^*$. □

For $\pi_{X,C}$ to be in $Ky\ \mu^*$, $X \in Ky\ \mu$, $X \ne C$, $X$ does not determine $C$ directly and $\alpha_{=C} \ne \emptyset$ for every cover $\alpha$.

### *Conclusion*

Many concepts from matroids and simplicial complexes have not been explored in this document. Representation theory and lattices were omitted. This is because the machinery of relational database theory required considerable space leaving little for other approaches. Representation needs to be explored further. There is a canonical relational representation that is like MSC representations ([R1, 4.4], [R2, 5.2]). These need to be compared.

We only gave a sufficient condition for membership of point sets like $\{(X, C)\}$ to be keys of $\mu^*$. A necessary condition is needed.

The remarkable observation of this paper is just how much structure of MSC carries over to general closures. Matroids are derive from algebraic systems that are extremely regular. General relational database dependencies are derived from the rules of physical and human systems. They are extremely irregular. Algebraic systems have computable dependencies with many symmetries. Relational dependencies are rarely computable and have few symmetries. However, the very fact that dependencies exist with closures imposes a consistent structure that is carried throughout all of these objects.

# *[R] References*